\let\WDDRefStepCounter\refstepcounter
\documentclass[final,oneeqnum]{siamltex}
\let\refstepcounter\WDDRefStepCounter

\usepackage{amsmath,amssymb}
\usepackage{mathrsfs}
\usepackage{mathtools}
\usepackage{bm}
\usepackage{graphicx}
\usepackage{booktabs}
\usepackage{cite}
\usepackage{hyperref}
\usepackage{cleveref}
\usepackage{aliascnt}
\usepackage{enumitem}
\usepackage{xcolor}
\usepackage{multirow}
\usepackage{placeins}
\usepackage{needspace}

\let\originalbibliography\bibliography
\renewcommand{\bibliography}[1]{\FloatBarrier\Needspace{8\baselineskip}{\raggedbottom\interlinepenalty=10000\originalbibliography{#1}}}

\hypersetup{hidelinks}
\graphicspath{{figures/}}


\newtheorem{definition}{Definition}[section]
\crefname{definition}{Definition}{Definitions}
\Crefname{definition}{Definition}{Definitions}

\newaliascnt{remark}{theorem}

\aliascntresetthe{remark}
\crefname{remark}{Remark}{Remarks}
\Crefname{remark}{Remark}{Remarks}

\newaliascnt{example}{theorem}

\aliascntresetthe{example}
\crefname{example}{Example}{Examples}
\Crefname{example}{Example}{Examples}

\newcommand{\C}{\mathbb{C}}
\newcommand{\N}{\mathbb{N}}

\newcommand{\bx}{\bm{x}}
\newcommand{\bz}{\bm{z}}
\newcommand{\bv}{\bm{v}}

\newcommand{\abs}[1]{\left|#1\right|}

\AddToHook{env/proof/before}{\Needspace{12\baselineskip}}
\AddToHook{env/definition/before}{\Needspace{14\baselineskip}}
\begin{document}

\title{
Efficient Numerical Evaluation Of Arbitrary High-order Partial Derivatives via Algebraic Decomposition
}

\author{Shuixin Fang\textsuperscript{\dag} \and Hang Zeng\textsuperscript{\dag}}
\maketitle
\begingroup
\renewcommand{\thefootnote}{\fnsymbol{footnote}}%
\footnotetext[2]{SKLMS \& Institute of Computational Mathematics and Scientific/Engineering Computing, Academy of Mathematics and Systems Science, Chinese Academy of Sciences, Beijing 100190, China.}%
\endgroup

\begin{abstract}
    Computing high-order partial derivatives of multivariate functions is often computationally and memory intensive.
    This paper addresses this problem via an algebraic approach.
    High-order partial derivatives of a multivariate function are expressed as linear combinations of univariate directional derivatives, which are computationally friendly.
    The coefficients and direction vectors are determined by the Waring decomposition of monomials, so that the number of univariate directional derivatives in the representation is minimized.
    Numerical experiments confirm the efficiency of the proposed method.
\end{abstract}



\pagestyle{myheadings}

\section{Introduction}\label{sec_intro}

Evaluating high-order partial derivatives of multivariate functions, such as deep neural networks, is a fundamental task in scientific computing.
As the derivative order and the input dimension increase, automatic differentiation~\cite{Griewank2008,Baydin2018,Paszke2019} becomes increasingly expensive in both computation and memory.
Univariate Taylor-mode automatic differentiation is a promising approach to address this challenge.
This approach recovers high-order partial derivatives as linear combinations of univariate directional derivatives, with suitably chosen coefficients and directions~\cite{Shi2024stde,Shi2026stdepp}.

Since each directional derivative still requires a pass of automatic differentiation, both computational cost and memory usage depend strongly on the number of directional derivatives in the representation.
A key question is therefore how to represent a given high-order partial derivative as a linear combination of as few univariate directional derivatives as possible.
This question motivates the present work.

Our approach is algebraic.
Treating the derivative operators as indeterminates of a polynomial ring, we reformulate the problem as expressing a homogeneous polynomial as a linear combination of as few powers of linear forms as possible.
This is precisely the Waring decomposition problem~\cite{CCG2012,BBT2013}, whose minimal decompositions yield representations of the given high-order partial derivative using the fewest univariate directional derivatives.
Numerical experiments show that the proposed Waring decomposition of derivatives (WDD) offers a significant efficiency advantage when the derivative contains many repeated indices.
We expect the approach to extend to partial derivatives of mixed order and to enable further study of minimal computational and memory costs.

\section{Directional representations of derivatives}\label{sec_formulation}

Let $\N_0$ denote the nonnegative integers and $\partial_j$ denote differentiation with respect to the input coordinate $x_j$. We write $\partial := (\partial_1, \ldots, \partial_d)^\top$.
For a multi-index $\alpha = (\alpha_1, \alpha_2, \ldots, \alpha_d) \in \N_0^d$, we define its order, factorial, and differential operator as
\begin{equation}\label{eq_multiindex}
    \abs{\alpha} := \sum_{j = 1}^{d} \alpha_j, \quad\quad \alpha! := \prod_{j = 1}^{d} \alpha_j!, \quad\quad \partial^{\alpha} := \partial_1^{\alpha_1} \partial_2^{\alpha_2} \cdots \partial_d^{\alpha_d}.
\end{equation}
For $p \in \N_0$, define the set of multi-indices of order $p$ as $\mathcal{M}_p := \{\alpha \in \N_0^d : \abs{\alpha} = p\}$.
For a direction $\bv = (v_1, v_2, \ldots, v_d)^\top \in \C^d$, set $\bv^\alpha := v_1^{\alpha_1} v_2^{\alpha_2} \cdots v_d^{\alpha_d}$.
For the zero multi-index $\boldsymbol{0}$, we write $\partial^{\boldsymbol{0}} u = u$.

Let $u : \C^d \to \C$ be holomorphic in a neighborhood of a given point $\bx \in \C^d$.
We seek to evaluate the derivative $\sum_{\alpha \in \mathcal{M}_p} A_{\alpha} \partial^\alpha u(\bx)$, where the positive integer $p$ and the coefficients $A_{\alpha} \in \C$ for $\alpha \in \mathcal{M}_p$ are given.
Assume that there exist a positive integer $R$, coefficients $c_i \in \C$, and directions $\bv_i \in \C^d$, $i = 1, 2, \ldots, R$, such that
\begin{equation}\label{eq_bridge}
    \sum_{\alpha \in \mathcal{M}_p} A_{\alpha} \bz^\alpha = \frac{1}{p!} \sum_{i = 1}^{R} c_i (\bv_i^\top \bz)^p \quad \text{for all } \bz \in \C^d.
\end{equation}
Since the constant-coefficient operators $\partial_j$ commute, substituting $z_j \mapsto \partial_j$ in \eqref{eq_bridge} gives
\begin{equation}\label{eq_schedule}
    \sum_{\alpha \in \mathcal{M}_p} A_{\alpha} \partial^\alpha u(\bx) = \frac{1}{p!} \sum_{i = 1}^{R} c_i (\bv_i^\top \partial)^p u(\bx)  = \frac{1}{p!} \sum_{i = 1}^{R} c_i \left.\frac{d^p}{d s^p} u(\bx + s \bv_i)\right|_{s = 0},
\end{equation}
where the second equality follows from the chain rule.
See \cite[Sec.~2]{CCG2012} and \cite[Secs.~3.7 and~3.17]{LoomisSternberg1990} for more details.

The right-hand side of \eqref{eq_schedule} can be evaluated by applying automatic differentiation to each univariate function $s \mapsto u(\bx + s \bv_i)$ \cite{Shi2024stde}.
When $R$ is small, this is expected to be more efficient than directly evaluating $\sum_{\alpha \in \mathcal{M}_p} A_{\alpha} \partial^\alpha u(\bx)$.
The remaining key issue is how to find coefficients $c_i$ and directions $\bv_i$, $i = 1, 2, \ldots, R$, that satisfy \eqref{eq_bridge} with a small $R$.
We give two constructions for the case of a single mixed partial derivative, i.e., $\sum_{\alpha \in \mathcal{M}_p} A_{\alpha} \partial^\alpha = \partial^\nu$ for some $\nu \in \mathcal{M}_p$.

The first construction follows from the inclusion--exclusion formula reviewed in \cite[Sec.~6.1, Eq.~(25)]{Shi2026stdepp}:
\begin{equation*}
    c_i = (-1)^{p - \abs{\gamma_i}} \prod_{j = 1}^{d} \binom{\nu_j}{\gamma_{i, j}}, \quad \bv_i = (\gamma_{i, 1}, \ldots, \gamma_{i, d})^\top, \quad i = 1, 2, \ldots, R,
\end{equation*}
where $\gamma_1, \ldots, \gamma_R$ enumerate the set
\begin{equation*}
    \bigl\{\gamma \in \N_0^d : \boldsymbol{0} \ne \gamma \leq \nu\bigr\}, \quad R = \prod_{j = 1}^{d} (\nu_j + 1) - 1.
\end{equation*}
Here $\gamma \leq \nu$ means $\gamma_j \leq \nu_j$ for all $j = 1, 2, \ldots, d$.
The resulting value of $R$ is not minimal for $p \geq 2$; see Corollary~\ref{thm_main} below.

The second construction is the Waring decomposition from \cite[Sec.~2, Eq.~(8)]{BBT2013}.
The pairs $(c_i, \bv_i)$, $i = 1, 2, \ldots, R$, enumerate the set $\{(c_{\bv}, \bv) : \bv \in \mathcal{V}\}$, where
\begin{equation*}
    c_{\bv} := \frac{\nu!}{R} \prod_{j \in I_1} v_j, \quad \mathcal{V} := \bigl\{\bv \in \C^d : v_j = 0,\, \forall j \in I_0,\quad v_j \in U_{\nu_j + 1},\, \forall j \in I_1,\quad v_{j_*} = 1\bigr\}.
\end{equation*}
The involving notations are defined by $I_0 := \{j \in \{1, \ldots, d\} : \nu_j = 0\}$, $I_1 := \{j \in \{1, \ldots, d\} : \nu_j > 0\}$, and 
\begin{align}
    &R = \frac{\prod_{j \in I_1} (\nu_j + 1)}{m_\nu + 1}, \quad j_* := \min T_\nu, \quad T_\nu := \{j \in I_1 : \nu_j = m_\nu\}, \label{eq_waring_weights} \\
    &m_\nu := \min_{j \in I_1} \nu_j, \quad U_m := \bigl\{\exp(2 \pi \mathrm{i} \ell / m) : \ell = 0, 1, \ldots, m - 1\bigr\},\quad m = 1, 2, \ldots. \label{eq_root_index_set}
\end{align}
Corollary~\ref{thm_main} below justifies the construction \eqref{eq_waring_weights}--\eqref{eq_root_index_set}; it is a direct consequence of the monomial Waring rank formula~\cite[Cor.~3.3 and Remark~2.3]{CCG2012}.

\begin{corollary}[Minimum representation length]\phantomsection\label{thm_main}
    Assume $A_{\nu} = 1$ and $A_{\alpha} = 0$ for $\alpha \in \mathcal{M}_p \setminus \{\nu\}$.
    The smallest positive integer $R$ for which there exist $c_i \in \C$ and $\bv_i \in \C^d$, $i = 1, \ldots, R$, satisfying \eqref{eq_bridge} is
    \begin{equation}\label{eq_directional_rank}
        R_\nu := \frac{\prod_{j \in I_1} (\nu_j + 1)}{m_\nu + 1}.
    \end{equation}
    The construction \eqref{eq_waring_weights}--\eqref{eq_root_index_set} attains this minimum.
\end{corollary}



\FloatBarrier
\section{Numerical experiments}\label{sec_experiments}

We compare WDD with a centered-polarization implementation of the directional Taylor identity reviewed in STDE++~\cite[Section~6.1, Eq.~(25)]{Shi2026stdepp} for individual partial derivatives $\partial_{i_1} \partial_{i_2} \cdots \partial_{i_p} u$.
We set the batch size to $B = 1000$.
The test network has 16 inputs and three hidden layers of width 128 with $\tanh$ activations, and both methods use identical weights and inputs and the same Taylor propagation routine on a Tesla V100.
For the real-valued networks and real inputs used here, we reuse Taylor coefficients at conjugate directions.

The order-$p$ coefficient counts in Table~\ref{tab_individual_t4} are per input point, after conjugate reuse in WDD and merging redundant directions in STDE++.
Times are the mean and sample standard deviation of three seed medians (30 synchronized calls after 10 warm-ups per seed), excluding direction preparation and parameter backpropagation.

\begin{table}[htbp]
    \centering
    \footnotesize
    \setlength{\tabcolsep}{2pt}
    \renewcommand{\arraystretch}{1.3}
    \caption{Ranks, Taylor coefficient counts, evaluation times, and peak allocated GPU memory on a Tesla V100 for WDD and STDE++~\cite[Section~6.1]{Shi2026stdepp}.
        Times are means and sample standard deviations of three seed medians, and speedups are means of paired ratios of STDE++ to WDD times.
        Memory values are maxima over three seeds in 30-s measurements after ten warm-up calls.}
    \label{tab_individual_t4}
    \resizebox{\linewidth}{!}{%
    \begin{tabular}{@{}lrrrrrrrrr@{}}
        \toprule
        Target & $p$ & $R$ & \multicolumn{2}{c}{Taylor coefficients} & \multicolumn{2}{c}{Time (ms)} & \multicolumn{2}{c}{Memory (MiB)} & Speedup \\
        \cmidrule(lr){4-5}\cmidrule(lr){6-7}\cmidrule(lr){8-9}
        & & & STDE++ & WDD & STDE++ & WDD & STDE++ & WDD & \\
        \midrule
        $\partial_1 \partial_2 \partial_3 \partial_4 \partial_5 \partial_6 u$ & 6 & 32 & 32 & 32 & $29.605 \pm 0.186$ & $29.600 \pm 0.145$ & 725.2 & 725.2 & $1.00 \times$ \\
        $\partial_1^{2} \partial_2^{2} \partial_3^{2} \partial_4^{2} u$ & 8 & 27 & 40 & 14 & $60.533 \pm 0.028$ & $45.456 \pm 0.029$ & 1172.9 & 834.4 & $1.33 \times$ \\
        $\partial_1^{3} \partial_2^{3} \partial_3^{3} u$ & 9 & 16 & 28 & 10 & $51.341 \pm 0.007$ & $39.932 \pm 0.019$ & 929.7 & 671.8 & $1.29 \times$ \\
        $\partial_1^{4} \partial_2^{4} \partial_3^{4} u$ & 12 & 25 & 49 & 13 & $153.037 \pm 0.025$ & $85.421 \pm 0.009$ & 2095.6 & 1129.5 & $1.79 \times$ \\
        $\partial_1^{5} \partial_2^{5} \partial_3^{5} u$ & 15 & 36 & 100 & 20 & $479.972 \pm 0.006$ & $195.631 \pm 0.026$ & 5260.4 & 2131.7 & $2.45 \times$ \\
        \bottomrule
    \end{tabular}%
    }
\end{table}

\FloatBarrier

\section{Conclusion and future work}\label{sec_concl}

We have represented a high-order partial derivative as a linear combination of univariate directional derivatives, with coefficients and directions given by a Waring decomposition of the corresponding monomial.
For a single mixed partial derivative $\partial^\nu$, the monomial Waring rank formula gives the minimum number $R_\nu$ of directional derivatives in \eqref{eq_schedule}, and the construction \eqref{eq_waring_weights}--\eqref{eq_root_index_set} attains it; see Corollary~\ref{thm_main}. 
In the experiments of Section~\ref{sec_experiments}, WDD matches the STDE++ baseline when all indices of $\partial^\nu$ are distinct and reduces the number of Taylor coefficients, the evaluation time, and the peak memory when indices are repeated.

Constructing minimal Waring decompositions for general homogeneous linear combinations of partial derivatives is a natural next step.
Let coefficients $A_\alpha$ be given for $\alpha \in \bigcup_{k=1}^p \mathcal{M}_k$.
Applying \eqref{eq_schedule} separately at each order yields a positive integer $R$, coefficients $c_i \in \C$, directions $\bv_i \in \C^d$, and orders $p_i \in \{1, \ldots, p\}$ for $i = 1, \ldots, R$, such that
\begin{equation*}
    \sum_{k = 1}^{p} \sum_{\alpha \in \mathcal{M}_k} A_{\alpha} \partial^\alpha u(\bx) = \sum_{i = 1}^{R} c_i \left.\Bigl(\frac{d}{d s}\Bigr)^{p_i} u(\bx + s \bv_i)\right|_{s = 0}.
\end{equation*}
Future work should study how to choose $\bv_i$, $c_i$, and $p_i$ so as to minimize the computational and memory costs of evaluating the right-hand side, including reusing Taylor coefficients across shared directions.

\FloatBarrier
\bibliographystyle{siamfull}
\bibliography{jsc_paper_main}

\end{document}